\documentclass[11pt,a4paper]{article}

\usepackage{inputenc}
\usepackage{amsmath}
\usepackage{bm}
\usepackage{bbold}
\usepackage{amsthm}

\usepackage{curves}
\usepackage{epic}

\usepackage{hyperref}

\title{Evaluation of the Lyapunov Exponent\\ for Generalized Linear Second-Order\\ Exponential Systems\thanks{Proc. 6th St.~Petersburg Workshop on Simulation, Volume II / Ed. by S.~M.~Ermakov and V.~B.~Melas and A.~N.~Pepelyshev, St.~Petersburg, 2009, pp.~875--881.}}

\author{Nikolai Krivulin\thanks{Faculty of Mathematics and Mechanics, St.~Petersburg State University, 28 Universitetsky Ave., St.~Petersburg, 198504, Russia, 
nkk@math.spbu.ru.} \thanks{The work was partially supported by the Russian Foundation for Basic Research under Grant \#09-01-00808.}
}

\date{}

\begin{document}

\maketitle

\begin{abstract}
We consider generalized linear stochastic dynamical systems with second-order state transition matrices. The entries of the matrix are assumed to be either independent and exponentially distributed or equal to zero. We give an overview of new results on evaluation of asymptotic growth rate of the system state vector, which is called the Lyapunov exponent of the system.
\end{abstract}

\section{Introduction}

The evolution of actual systems that occur in management, engineering, computer sciences, and other areas can frequently be represented through stochastic dynamic equations of the form
$$
\bm{z}(k)=A(k)\bm{z}(k-1),
$$
where $ A(k) $ is a random state transition matrix, $ \bm{z}(k) $ is a system state vector, and matrix-vector multiplication is thought of as defined in terms of a semiring with the operations of taking maximum and addition \cite{Kolokoltsov1997Idempotent,Litvinov1998Idempotent,Heidergott2006Max-plus}.

In many cases, the analysis of a system involves evaluation of asymptotic growth rate of the system state vector $ \bm{z}(k) $, which is normally referred to as the Lyapunov exponent \cite{Heidergott2007Max-plus,Jean-Marie1994Analytical}.

Evaluation of the Lyapunov exponent typically appears to be a difficult problem even for quite simple systems. Related results include the solutions obtained in \cite{Olsder1990Discrete,Jean-Marie1994Analytical} for systems with matrices of the second order with independent and exponentially distributed entries. In \cite{Olsder1990Discrete}, the Lyapunov exponent is obtained in the case that all entries of the matrix are identically distributed with unit mean.

Further results are given in \cite{Jean-Marie1994Analytical} under the condition that the diagonal entries have one common distribution, whereas the off-diagonal entries do follow another distribution. A system with a matrix such that its diagonal entries are distributed with unit mean, and the off-diagonal entries are equal to zero is also examined.

The purpose of this paper is to give an overview of new results which are related to evaluation of the Lyapunov exponent in generalized linear systems that have matrices of the second order with exponentially distributed entries (second-order exponential systems).

\section{Stochastic Linear Dynamical System}

Consider a dynamical system that can be represented through the linear equation in the semiring with the operations of maximum and addition
$$
\bm{z}(k)=A(k)\bm{z}(k-1),
$$
where
$$
A(k)
=
\left(
	\begin{array}{cc}
		\alpha_{k}	& \beta_{k} \\
		\gamma_{k}	& \delta_{k}
	\end{array}
\right),
\qquad
\bm{z}(k)
=
\left(
	\begin{array}{c}
		x(k) \\
		y(k)
	\end{array}
\right),
\qquad
\bm{z}(0)
=
\left(
	\begin{array}{c}
		0 \\
		0
	\end{array}
\right).
$$

With ordinary notation, the above vector equation can be written as
\begin{align*}
x(k)&=\max(x(k-1)+\alpha_{k},y(k-1)+\beta_{k}), \\
y(k)&=\max(x(k-1)+\gamma_{k},y(k-1)+\delta_{k}).
\end{align*}

The Lyapunov exponent for the system is given by
$$
\lambda=\lim_{k\to\infty}\frac{1}{k}\max(x(k),y(k)).
$$

Suppose the sequences $ \{\alpha_{k} \} $, $ \{\beta_{k} \} $, $ \{\gamma_{k} \} $, and $ \{\delta_{k}\} $ each involve independent and identically distributed random variables; $ \alpha_{k} $, $ \beta_{l} $, $ \gamma_{m} $, and $ \delta_{n} $ are independent for any $ k,l,m,n $. Finally, we assume that $ \alpha_{k} $, $ \beta_{k} $, $ \gamma_{k} $, and $ \delta_{k} $ have the exponential probability distributions with respective parameters $ \mu $, $ \nu $, $ \sigma $, and $ \tau $.

\section{Systems With Matrices Having Zero Entries}

We start with results obtained in \cite{Krivulin2007Growth,Krivulin2009Calculating} for systems with state transition matrices having one or two nonrandom entries that are equal to zero. The reduced number of random entries in the matrices allows one to simplify the evaluation of the Lyapunov exponent and usually gives quite compact results.

The solution method is based on construction of a sequence of probability distribution functions. The convergence of the sequence is examined and the limiting distribution is derived as the solution of an integral equation. The Lyapunov exponent is then evaluated as the expected value of a random variable determined through the limiting distribution function.

\subsection{Matrix With Zero Off-Diagonal Entries}

The system state transition matrix together with its related result take the form 
$$
A(k)
=
\left(
	\begin{array}{cc}
		\alpha_{k}	& 0 \\
		0						& \delta_{k}
	\end{array}
\right),
\qquad
\lambda
=
\frac{\mu^{4}+\mu^{3}\tau+\mu^{2}\tau^{2}+\mu\tau^{3}+\tau^{4}}{\mu\tau(\mu+\tau)(\mu^{2}+\tau^{2})}.
$$

With $ \tau=\mu=1 $ we have the result $ \lambda=1.25 $ which coincides with that in \cite{Jean-Marie1994Analytical}.

\subsection{Matrix With Zero Diagonal}

In this case, the matrix and the Lyapunov exponent are represented as
$$
A(k)
=
\left(
	\begin{array}{cc}
		0						& \beta_{k} \\
		\gamma_{k}	& 0
	\end{array}
\right),
\qquad
\lambda
=
\frac{4\nu^{2}+7\nu\sigma+4\sigma^{2}}{6\nu\sigma(\nu+\sigma)}.
$$

\subsection{Matrix With Zero Row or Column}

Provided that the second row in the matrix has only zero entries, we arrive at
$$
A(k)
=
\left(
	\begin{array}{cc}
		\alpha_{k}	& \beta_{k} \\
		0						& 0
	\end{array}
\right),
\qquad
\lambda
=
\frac{2\mu^{4}+7\mu^{3}\nu+10\mu^{2}\nu^{2}+11\mu\nu^{3}+4\nu^{4}}{\mu\nu(\mu+\nu)^{2}(3\mu+4\nu)}.
$$

When the entries of the second column are zero, we have
$$
A(k)
=
\left(
	\begin{array}{cc}
		\alpha_{k}	& 0 \\
		\gamma_{k}	& 0
	\end{array}
\right),
\qquad
\lambda
=
\frac{2\mu^{4}+7\mu^{3}\sigma+10\mu^{2}\sigma^{2}+11\mu\sigma^{3}+4\sigma^{4}}{\mu\sigma(\mu+\sigma)^{2}(3\mu+4\sigma)}.
$$

\subsection{Matrix With Zero Entry on Diagonal}

Consider a system with the state transition matrix
$$
A(k)
=
\left(
	\begin{array}{cc}
		\alpha_{k}	& \beta_{k} \\
		\gamma_{k}	& 0
	\end{array}
\right).
$$

Whereas evaluation of the Lyapunov exponent for this system in the general case leads to rather cumbersome algebraic manipulations, there are two main particular cases which offer their related results in a relatively compact form. Under the condition that $ \sigma=\mu $, we have
$$
\lambda
=
\frac{48\mu^{5}+238\mu^{4}\nu+495\mu^{3}\nu^{2}+581\mu^{2}\nu^{3}+326\mu\nu^{4}+68\nu^{5}}{2\mu\nu(36\mu^{4}+147\mu^{3}\nu+215\mu^{2}\nu^{2}+130\mu\nu^{3}+28\nu^{4})}.
$$

Provided that $ \sigma=\nu $, the value of the Lyapunov exponent is given by
$$
\lambda=P(\mu,\nu)/Q(\mu,\nu),
$$
where
\begin{align*}
P(\mu,\nu)
&=
15\mu^{8}+152\mu^{7}\nu+624\mu^{6}\nu^{2}+1382\mu^{5}\nu^{3}+1838\mu^{4}\nu^{4}+1592\mu^{3}\nu^{5}
\\
&\qquad+
973\mu^{2}\nu^{6}+384\mu\nu^{7}+64\nu^{8},
\\
Q(\mu,\nu)
&=
\mu\nu(\mu+\nu)^{2}(12\mu^{5}+97\mu^{4}\nu+286\mu^{3}\nu^{2}+397\mu^{2}\nu^{3}+256\mu\nu^{4}+64\nu^{5}).
\end{align*}

\subsection{Matrix With Zero Entry Below Diagonal}

Suppose that there is a system with the state transition matrix defined as
$$
A(k)
=
\left(
	\begin{array}{cc}
		\alpha_{k}	& \beta_{k} \\
		0						& \delta_{k}
	\end{array}
\right).
$$

Consider two particular cases. With the condition $ \nu=\mu $, we have the result
$$
\lambda
=
P(\mu,\tau)/Q(\mu,\tau),
$$
where
\begin{align*}
P(\mu,\tau)
&=
288\mu^{8}+1048\mu^{7}\tau+1936\mu^{6}\tau^{2}+2688\mu^{5}\tau^{3}+3012\mu^{4}\tau^{4}
\\
&\qquad
+2226\mu^{3}\tau^{5}+941\mu^{2}\tau^{6}+204\mu\tau^{7}+17\tau^{8},
\\
Q(\mu,\tau)
&=
2\mu\tau(144\mu^{7}+524\mu^{6}\tau+968\mu^{5}\tau^{2}+1200\mu^{4}\tau^{3}+910\mu^{3}\tau^{4}
\\
&\qquad
+387\mu^{2}\tau^{5}+84\mu\tau^{6}+7\tau^{7}).
\end{align*}

Provided that $ \tau=\mu $, the solution takes the form
$$
\lambda
=
P(\mu,\nu)/Q(\mu,\nu),
$$
where
\begin{align*}
P(\mu,\nu)
&=
256\mu^{10}+2112\mu^{9}\nu+8044\mu^{8}\nu^{2}+19355\mu^{7}\nu^{3}+32167\mu^{6}\nu^{4}
\\
&\qquad
+36887\mu^{5}\nu^{5}+28709\mu^{4}\nu^{6}+14854\mu^{3}\nu^{7}+4912\mu^{2}\nu^{8}
\\
&\qquad
+4\mu\nu^{9}+80\nu^{10},
\\
Q(\mu,\nu)
&=
2\mu\nu(\mu+\nu)(192\mu^{8}+1344\mu^{7}\nu+4047\mu^{6}\nu^{2}+6770\mu^{5}\nu^{3}
\\
&\qquad
+6799\mu^{4}\nu^{4}+4216\mu^{3}\nu^{5}+1600\mu^{2}\nu^{6}+344\mu\nu^{7}+32\nu^{8}).
\end{align*}

\section{General Second-Order Exponential System}

Consider a general second-order exponential system which has the matrix
$$
A(k)
=
\left(
\begin{array}{cc}
	\alpha_{k}	& \beta_{k} \\
	\gamma_{k}	& \delta_{k}
\end{array}
\right)
$$
with its entries $ \alpha_{k} $, $ \beta_{k} $, $ \gamma_{k} $, and $ \delta_{k} $ assumed to be independent random variables that are exponentially distributed with the respective parameters $ \mu $, $ \nu $, $ \sigma $, and $ \tau $.

To get the Lyapunov exponent, a computational technique developed in \cite{Krivulin2008Evaluation} can be implemented which reduces the problem to the solution of a system of linear equations, accompanied by the evaluation of a linear functional of the system solution. The technique leans upon construction and examination of a sequence of probability density functions. It is shown that there is a one-to-one correspondence between the density functions and vectors in a vector space. The correspondence is then exploited to provide for the solution in terms of algebraic computations. 

Based on the above technique, the Lyapunov exponent can be evaluated as follows. First we introduce the vectors
$$
\bm{\omega}_{1}
=
(\omega_{10},\omega_{11},\omega_{12},\omega_{13})^{T},
\quad
\bm{\omega}_{2}
=
(\omega_{20},\omega_{21},\omega_{22},\omega_{23})^{T},
\quad
\bm{\omega}
=
(\bm{\omega}_{1}^{T},\bm{\omega}_{2}^{T})^{T}.
$$

Furthermore, we define the matrices
$$
U_{1}
=
\left(
\begin{array}{ccc}
1 & 1 & 1 \\
\frac{\mu}{\mu+\nu} & \frac{1}{2} & \frac{\mu+\nu}{\mu+2\nu} \\
\frac{\mu}{\mu+\tau} & \frac{\nu}{\nu+\tau} & \frac{\mu+\nu}{\mu+\nu+\tau} \\
\frac{\mu}{\mu+\nu+\tau} & \frac{\nu}{2\nu+\tau} & \frac{\mu+\nu}{\mu+2\nu+\tau}
\end{array}
\right),
\qquad
U_{2}
=
\left(
\begin{array}{ccc}
1 & 1 & 1 \\
\frac{\sigma}{\mu+\sigma} & \frac{\tau}{\mu+\tau} & \frac{\sigma+\tau}{\mu+\sigma+\tau} \\
\frac{1}{2} & \frac{\tau}{\sigma+\tau} & \frac{\sigma+\tau}{2\sigma+\tau} \\
\frac{\sigma}{\mu+2\sigma} & \frac{\tau}{\mu+\sigma+\tau} & \frac{\sigma+\tau}{\mu+2\sigma+\tau}
\end{array}
\right),
$$
\begin{align*}
V_{11}
&=
\left(
\begin{array}{cccc}
\frac{\sigma}{\mu+\sigma} & 0 & -\frac{\mu\sigma}{(\mu+\tau)(\mu+\sigma+\tau)} & 0 \\
0 & \frac{\sigma}{\nu+\sigma} & 0 & -\frac{\nu\sigma}{(\nu+\tau)(\nu+\sigma+\tau)} \\
0 & -\frac{\sigma}{\mu+\nu+\sigma} & 0 & \frac{\sigma(\mu+\nu)}{(\mu+\nu+\tau)(\mu+\nu+\sigma+\tau)}
\end{array}
\right),
\\
V_{12}
&=
\left(
\begin{array}{cccc}
0 & \frac{\tau}{\mu+\tau} & 0 & -\frac{\mu\tau}{(\mu+\sigma)(\mu+\sigma+\tau)} \\
\frac{\tau}{\nu+\tau} & 0 & -\frac{\nu\tau}{(\nu+\sigma)(\nu+\sigma+\tau)} & 0 \\
0 & -\frac{\tau}{\mu+\nu+\tau} & 0 & \frac{\tau(\mu+\nu)}{(\mu+\nu+\sigma)(\mu+\nu+\sigma+\tau)}
\end{array}
\right),
\\
V_{21}
&=
\left(
\begin{array}{cccc}
\frac{\mu}{\mu+\sigma} & -\frac{\mu\sigma}{(\nu+\sigma)(\mu+\nu+\sigma)} & 0 & 0 \\
0 & 0 & \frac{\mu}{\mu+\tau} & -\frac{\mu\tau}{(\nu+\tau)(\mu+\nu+\tau)} \\
0 & 0 & -\frac{\mu}{\mu+\sigma+\tau} & \frac{\mu(\sigma+\tau)}{(\nu+\sigma+\tau)(\mu+\nu+\sigma+\tau)}
\end{array}
\right),
\\
V_{22}
&=
\left(
\begin{array}{cccc}
0 & 0 & \frac{\nu}{\nu+\sigma} & -\frac{\nu\sigma}{(\mu+\sigma)(\mu+\nu+\sigma)} \\
\frac{\nu}{\nu+\tau} & -\frac{\nu\tau}{(\mu+\tau)(\mu+\nu+\tau)} & 0 & 0 \\
0 & 0 & -\frac{\nu}{\nu+\sigma+\tau} & \frac{\nu(\sigma+\tau)}{(\mu+\sigma+\tau)(\mu+\nu+\sigma+\tau)}
\end{array}
\right).
\end{align*}

Suppose that the vector $ \bm{\omega} $ is the solution of the system
\begin{align*}
(I-W)\bm{\omega}
&=
\bm{0}, \\
\omega_{10}+\omega_{20}
&=
1,
\end{align*}
where $ I $ represents identity matrix, and
$$
W
=
\left(
\begin{array}{cc}
U_{1}V_{11} & U_{1}V_{12} \\
U_{2}V_{21} & U_{2}V_{22}
\end{array}
\right).
$$

The value of the Lyapunov exponent is then given by
$$
\lambda
=
\bm{q}_{1}^{T}\bm{\omega}_{1}+\bm{q}_{2}^{T}\bm{\omega}_{2},
$$
where $ \bm{q}_{1} $ and $ \bm{q}_{2} $ are vectors such that
\begin{align*}
\bm{q}_{1}
&=
\left(
\begin{array}{cc}
\frac{\mu^{2}+\mu\sigma+\sigma^{2}}{\mu\sigma(\mu+\sigma)} \\
\frac{\mu\sigma(\mu+2\nu+\sigma)}{\nu(\mu+\nu)(\nu+\sigma)(\mu+\nu+\sigma)} \\
\frac{\mu\sigma(\mu+\sigma+2\tau)}{\tau(\mu+\tau)(\sigma+\tau)(\mu+\sigma+\tau)} \\
-
\frac{\mu\sigma(\mu+2\nu+2\tau+\sigma)}{(\nu+\tau)(\mu+\nu+\tau)(\nu+\sigma+\tau)(\mu+\nu+\sigma+\tau)}
\end{array}
\right),
\\
\bm{q}_{2}
&=
\left(
\begin{array}{c}
\frac{\nu^{2}+\nu\tau+\tau^{2}}{\nu\tau(\nu+\tau)} \\
\frac{\nu\tau(2\mu+\nu+\tau)}{\mu(\mu+\nu)(\mu+\tau)(\mu+\nu+\tau)} \\
\frac{\nu\tau(\nu+2\sigma+\tau)}{\sigma(\nu+\sigma)(\sigma+\tau)(\nu+\sigma+\tau)} \\
-
\frac{\nu\tau(2\mu+\nu+2\sigma+\tau)}{(\mu+\sigma)(\mu+\nu+\sigma)(\mu+\sigma+\tau)(\mu+\nu+\sigma+\tau)}
\end{array}
\right).
\end{align*}

Suppose that $ \tau=\mu $ and $ \sigma=\nu $. Implementation of the above technique gives the solution
$$
\lambda=P(\mu,\nu)/Q(\mu,\nu),
$$
where
\begin{align*}
P(\mu,\nu)
&=
160\mu^{10}+1776\mu^{9}\nu+8220\mu^{8}\nu^{2}+21378\mu^{7}\nu^{3}+35595\mu^{6}\nu^{4}
\\
&\qquad+
41566\mu^{5}\nu^{5}+35595\mu^{4}\nu^{6}+21378\mu^{3}\nu^{7}+8220\mu^{2}\nu^{8}
\\
&\qquad+
1776\mu\nu^{9}+160\nu^{10},
\\
Q(\mu,\nu)
&=
16\mu\nu(\mu+\nu)(8\mu^{8}+80\mu^{7}\nu+321\mu^{6}\nu^{2}+690\mu^{5}\nu^{3}+880\mu^{4}\nu^{4}
\\
&\qquad+
690\mu^{3}\nu^{5}+321\mu^{2}\nu^{6}+80\mu\nu^{7}+8\nu^{8}).
\end{align*}

Note that the obtained solution coincides with that in \cite{Jean-Marie1994Analytical}.

The author is grateful to the anonymous reviewer for valuable comments and suggestions.

\bibliographystyle{utphys}

\bibliography{Evaluation_of_the_Lyapunov_exponent_for_generalized_linear_second-order_exponential_systems}

\end{document}